\documentstyle[12pt]{article}
\title{Solutions of DEs and PDEs as Potential Maps Using First Order
Lagrangians}
\author{Constantin Udri\c ste \\ {} \\
University Politehnica of Bucharest \\
Department of Mathematics I \\
Splaiul Independen\c tei 313 \\
77206 Bucharest, Romania\\
email:udriste@mathem.pub.ro}

\def\interior{{\setlength{\unitlength}{1em}\begin{picture}(0.85,1)
	\put(.1,0){\line(1,0){.55}}\put(.65,0){\line(0,1){.55}}\end{picture}}}

\date{}
\begin{document}
\maketitle

\newcommand{\ty}{\infty}
\newcommand{\te}{\theta}
\newcommand{\de}{\delta}
\newcommand{\ov}{\over}
\newcommand{\di}{\displaystyle}
\newcommand{\si}{\sigma}
\newcommand{\na}{\nabla}
\newcommand{\pa}{\partial}
\newcommand{\la}{\lambda}
\newcommand{\al}{\alpha}
\newcommand{\ep}{\varepsilon}
\newcommand{\va}{\varphi}
\newcommand{\ld}{\ldots}
\newcommand{\noa}{\noalign{\medskip}}
\newcommand{\ti}{\times}

\begin{abstract}
Using parametrized curves (Section 1) or parametrized sheets (Section 3), and
suitable metrics, we treat the jet bundle of order one as a semi-Riemann
manifold. This point of view allows the description of solutions of DEs as
pregeodesics (Section 1) and the solutions of PDEs as potential maps (Section 3),
via Lagrangians of order one or via generalized Lorentz world-force laws.
Implicitly, we solved a problem rised first by Poincar\'e: find a suitable
geometric structure that converts the trajectories of a given vector field
into geodesics (see also [6] - [11]).
Section 2 and Section 3 realize the passage from the Lagrangian dynamics to
the covariant Hamilton equations.
\end{abstract}



\section{Solutions of DEs as pregeodesics}

Unless specifically denied, all manifolds, all objects on them, and all
maps from one manifold into another will be $C^\ty$; however, we sometimes
redundantly write "a $C^\ty$ manifold", and so on, for emphasis.

Let $(T = R,h)$ and $(M,g)$ be semi-Riemann manifolds of dimensions 1 and $n$.
Hereafter we shall assume that the manifold $T$ is oriented. Latin letters
will be used for indexing the components of geometrical objects attached to
the manifold $M$.

Local coordinates will be written
$$
t = t^1, \quad x = (x^i), \quad i = 1,\ld, n,
$$
and the components of the corresponding metric tensors and Christoffel
symbols will be denoted by $h_{11}$, $g_{ij}$, $H^1_{11}$, $G^i_{jk}$.
Indices of distinguished objects will be rised and lowered in the usual
fashion.

Let $C^\ty (T,M) = \{\va: T \to M\:|\: \va \; \hbox{of class}\; C^\ty\}$.
For any $\va,\psi \in C^\ty (T,M)$, we define the equivalence relation
$\va \sim \psi$ at $(t_0, x_0) \in T \times M$ by
$$
x^i (t_0) = y^i (t_0) = x^i_0, \quad {dx^i \ov dt} (t_0) = {dy^i \ov dt}
(t_0).
$$
Using the factorization
$$
J^1_{(t_0,x_0)} (T,M) = C^\ty (T,M)/\sim
$$
we introduce the jet bundle of order one
$$
J^1 (T,M) = \bigcup_{(t_0,x_0)\in T\times M} J^1_{(t_0,x_0)} (T,M).
$$
Denoting by $[\va]_{(t_0,x_0)}$ the equivalence class of the map $\va$, we
define the projection
$$
\pi: J^1 (T,M) \to T \times M, \quad
\pi[\va]_{(t_0,x_0)} = (t_0, \va (t_0)).
$$

Suppose that the base $T \times M$ is covered by a system of coordinate
neighborhoods $(U\times V, t^\al, x^i)$. Then we can define the
diffeomorphism
$$
F_{U\ti V}: \pi^{-1} (U \ti V) \to U \ti V \ti R^{1\cdot n}
$$
$$
F_{UV} [\va]_{(t_0,x_0)} = \left( t_0, x^i_0, {dx^i \ov dt} (t_0) \right).
$$

Consequently $J^1(T,M)$ is a differentiable manifold of dimension
$1+n+ 1\cdot n = 2n+1$. The coordinates on $\pi^{-1} (U \ti V) \subset
J^1 (T,M)$ will be
$$
\left( t^1 = t, x^i, y^i = {dx^i \ov dt}\right),
$$
where
$$
t^1 ([\va]_{(t_0,x_0)}) = t^1(t_0), x^i([\va]_{(t_0,x_0)}) = x^i(x_0), \;
y^i([\va]_{(t_0,x_0)}) = {dx^i \ov dt} (t_0).
$$
A local changing of coordinates $(t, x^i, y^i) \to (\bar t, \bar x^i,
\bar y^i)$ is given by
$$
\bar t = \bar t(t), \; \bar x^i = \bar x^i(x^j), \;
\bar y^i = {\pa \bar x^i \ov \pa x^j}{dt \ov d\bar t} \; y^j, \leqno (1)
$$
where
$$
{d\bar t \ov dt} > 0, \quad \det\left( {\pa \bar x^i \ov \pa x^j} \right)
\ne 0.
$$
The expression of the Jacobian matrix of the local diffeomorphism $(1)$ shows
that the jet bundle of order one $J^1 (T,M)$ is always orientable.

Let $H^1_{11} = \di{1 \ov 2} h^{11} \di{dh_{11} \ov dt^1} = \di{1 \ov 2}
h^{-1}_{11} \di{dh_{11} \ov dt^1} = \di{1 \ov 2} \di{d \ov dt^1}
\sqrt{|h_{11}|}$, $G^i_{jk}$ be the components of the connections induced
by $h$ and $g$ respectively. If $\di\left( t=t^1, x^i, y^i = {dx^i \ov dt}\right)$
are the coordinates of a point in $J^1(T,M)$, then
$$
{\delta \ov dt}{dx^i \ov dt} =
{d^2x^i \ov dt^2} - H^1_{11} {dx^i \ov dt} + G^i_{jk} {dx^j \ov dt}
{dx^k \ov dt}
$$
are the components of a distinguished tensor on $T \times M$. Also
$$
\left( {\delta \ov \delta t} = {d \ov dt} + H^1_{11} y^i {\pa \ov \pa y^i},
\quad {\delta \ov \delta x^i} = {\pa \ov \pa x^i} - G^h_{ik} y^k
{\pa \ov \pa y^h}, \; {\pa \ov \pa y^i} \right)
$$
$$
(dt, \; dx^j, \; \delta y^j = dy^j - H^1_{11} y^j dt + G^j_{hk} y^h dx^k)
$$
are dual frames on $J^1(T,M)$, i.e.,
$$
dt\left( {\delta \ov \delta t} \right) = 1, \;
dt\left( {\delta \ov \delta x^i} \right) = 0, \;
dt\left( {\pa \ov \pa y^i} \right) = 0
$$
$$
dx^j\left( {\delta \ov \delta t} \right) = 0, \;
dx^j\left( {\delta \ov \delta x^i} \right) = \delta^j_i, \;
dx^j\left( {\pa \ov \pa y^i} \right) = 0
$$
$$
\delta y^j\left( {\delta \ov \delta t} \right) = 0, \;
\delta y^j\left( {\delta \ov \delta x^i} \right) = 0,\;
\delta y^j\left( {\pa \ov \pa y^i} \right) = \delta^j_i.
$$

Using these frames, we define on $J^1 (T,M)$ the induced Sasaki-like metric
$$
S_1 = h_{11} dt \otimes dt + g_{ij} dx^i \otimes dx^j + h^{11} g_{ij}
\delta y^i \otimes \delta y^j.
$$
The geometry of the manifold $(J^1(T,M), S_1)$ was developed recently in [4].

Now we shall generalize the Lorentz world-force law which was initially
stated [5] for particles in nonquantum relativity.

{\bf Definition}. Let $F = (F_j{}^i)$ and $U = (U^i)$ be $C^\ty$ distinguished
tensors on $T \times M$, where $\omega_{ji} = g_{hi} F_j{}^h$ is skew-symmetric
with respect to $j$ and $i$. Let $c(t,x)$ be a $C^\ty$ real function on
$T \times M$. A map $\va: T \to M$ obeys the {\it Lorentz-Udri\c ste World-Force
Law} with respect to $F, U, c$ iff
$$
h^{11} {\delta \ov dt} {dx^i \ov dt} = h^{11}
\left( g^{ij} {\pa c \ov \pa x^j} + F_j{}^i \; {dx^j \ov dt} + U^i \right).
$$

Now we remark that a $C^\ty$ distinguished tensor field $X^i(t,x)$,
$i=1,\ld, n$ on $T \ti M$ defines a family of trajectories as solutions of
DEs system of order one
$$
{dx^i \ov dt} = X^i (t, x(t)). \leqno (2)
$$
The distinguished tensor field $X^i(t,x)$ and semi-Riemann metrics $h$ and
$g$ determine the {\it potential energy}
$$
f: T \ti M \to R, \quad f = {1 \ov 2} h^{11} g_{ij} X^i X^j.
$$
The distinguished tensor field (family of trajectories) $X^i$ on
$(T\ti M, h_{11} + g)$ is called:

1) {\it timelike}, if $f < 0$;

2) {\it nonspacelike} or {\it causal}, if $f \le 0$;

3) {\it null} or {\it lightlike}, if $f=0$;

4) {\it spacelike}, if $f > 0$.

Let $X^i$ be a distinguished tensor field of everywhere constant energy.
If $X^i$ (the system (2)) has no critical point on $M$, then upon rescaling,
it may be supposed that $f \in \{-1,0,1\}$. Generally,
$\cal E = \{x_0 \in M| X^i(t,x_0) = 0, \forall t \in T\}$ is the set of critical
points of the distinguished tensor field, and this rescaling is possible
only on $T \ti (M\setminus \cal E)$.

Using the operator (derivative along a solution of (2))
$$
{\delta \ov dt} {dx^i \ov dt} = {d^2x^i \ov dt^2} -
H^1_{11} {dx^i \ov dt} + G^i_{jk} {dx^j \ov dt}{dx^k \ov dt},
$$
the Levy-Civita connection $D$ of $(R,h)$ and the Levy-Civita connection
$\nabla$ of $(M,g)$,  we obtain the prolongation (system of DEs of order two)
$$
{d^2x^i \ov dt^2} - H^1_{11} {dx^i \ov dt} + G^i_{jk} {dx^j \ov dt}{dx^k \ov dt}
= DX^i + (\nabla_j X^i){dx^j \ov dt}, \leqno (3)
$$
where
$$
\na_j X^i = {\pa X^i \ov\pa x^j} + G^i_{jk} X^k, \quad
DX^i = {\pa X^i \ov \pa t} - H^1_{11} X^i.
$$
The distinguished tensor field $X^i$, the metric $g$, and the connection
$\na$ determine the external distinguished tensor field
$$
F_j{}^i = \na_j X^i - g^{ih} g_{kj} \na_h X^k,
$$
which characterizes the {\it helicity} of the distinguished tensor field $X^i$.

The DEs system (3) can be written in the equivalent form
$$
{d^2x^i \ov dt^2} - H^1_{11} {dx^i \ov dt} + G^i_{jk} {dx^j \ov dt}{dx^k \ov dt}
= g^{ih} g_{kj} (\na_h X^k) {dx^j \ov dt} + F_j{}^i {dx^j \ov dt} + DX^i.
\leqno (4)
$$
Now we modify this DEs system into
$$
{d^2x^i \ov dt^2} - H^1_{11} {dx^i \ov dt} + G^i_{jk} {dx^j \ov dt}{dx^k \ov dt}
= g^{ih} g_{kj} (\na_h X^k)X^j + F_j{}^i {dx^j \ov dt} + DX^i
\leqno (5)
$$
or equivalently,
$$
h^{11} \left( {d^2x^i \ov dt^2} - H^1_{11} {dx^i \ov dt} + G^i_{jk} {dx^j \ov dt}{dx^k \ov dt}\right)
= h^{11} (g^{ih} g_{kj} (\na_h X^k) X^j + F_j{}^i {dx^j \ov dt} + DX^i).
$$
The system (5) is still a prolongation of the DEs system (2).

{\bf Theorem}. {\it The kinematic system (2) can be prolonged to the second
order dynamical system (5).}

{\bf Corollary}. {\it Choosing the metrics $h$ and $g$ such that $f \in \{-1,0,1\}$,
then the kinematic system (2) can be prolonged to the second order dynamical
system}
$$
{d^2x^i \ov dt^2} - H^1_{11} {dx^i \ov dt} + G^i_{jk} {dx^j \ov dt}{dx^k \ov dt}
= F_j{}^i {dx^j \ov dt} + DX^i.
$$

We shall show that the dynamical system (5) is in fact an Euler-Lagrange
system. We identify $J^1(T\ti M)$ with its dual via the semi-Riemann metrics
$h$ and $g$.

{\bf Theorem}. 1) {\it The solutions of the DEs system (5) are the extremals
of the Lagrangian}
$$
L = {1 \ov 2} h^{11} g_{ij} \left( {dx^i \ov dt} - X^i \right)
\left( {dx^j \ov dt} - X^j \right) \sqrt{|h_{11}|} =
$$
$$
= \left( {1 \ov 2} h^{11} g_{ij} {dx^i \ov dt} {dx^j \ov dt} - h^{11}g_{ij}
{dx^i \ov dt} X^j + f\right) \sqrt{|h_{11}|}.
$$

2) {\it If $F_j{}^i = 0$, then the solutions of the DEs system (5) are the
extremals of the Lagrangian}
$$
L = \left( {1 \ov 2} h^{11} g_{ij} {dx^i \ov dt} {dx^j \ov dt} + f\right)
\sqrt{|h_{11}|}.
$$

3) {\it Both Lagrangians produce the same Hamiltonian}
$$
H = \left( {1 \ov 2} h^{11} g_{ij} {dx^i \ov dt} {dx^j \ov dt} - f\right)
\sqrt{|h_{11}|}.
$$

{\bf Theorem (Lorentz-Udri\c ste World-Force Law)}.

1) {\it Every solution of DEs system
$$
{d^2x^i \ov dt^2} - H^1_{11} {dx^i \ov dt} + G^i_{jk} {dx^j \ov dt}
{dx^k \ov dt} = g^{ih}g_{kj} (\na_h X^k)X^j + DX^i
$$
is a pregeodesic (potential map) on the semi-Riemann manifold}
$(T\ti M, h+g)$.

2) {\it Every solution of DEs system (5) is a horizontal pregeodesic
(potential map) on the semi-Riemann-Lagrange manifold}
$$
(T\ti M, h+g, \;
N(^i_1)_j = G^i_{jk}y^k - F_j{}^i, \quad M(^i_1)_1 = - H^1_{11} y^i).
$$

{\bf Corollary}. {\it Every DE generates a Lagrangian of order one  via the
associated first order DEs system and suitable metrics on the manifold of
independent variable and on the manifold of functions. In this sense the
solutions of the initial DE are pregeodesics produced by a suitable
Lagrangian}.

{\bf Proof}. Let $t \in R$ denote a real variable, usually referred to as
the time. It may be pointed out that the DE
$$
{d^n x \ov dt^n} = f\left( t, x, {dx \ov dt}, \ldots, {d^{n-1} x \ov dt^{n-1}}
\right), \leqno (6)
$$
where $x$ is the unknown function, is equivalent to a system (2). For if we
set $x=x^1$, then (6) is equivalent to
$$
{dx^1 \ov dt} = x^2, \quad {dx^2 \ov dt} = x^3, \ld, {dx^{n-1} \ov dt} = x^n
$$
$$
{dx^{n} \ov dt} = f(t, x^1, x^2, \ld, x^n),
$$
which is type (2). Therefore, the preceding theory applies.

\section{Hamiltonian approach}

Let $(Q, \Omega)$ be a symplectic manifold (of even dimension). The Hamiltonian
vector field $X_H$ of the function $H \in \cal F(Q)$ is defined by
$$
X_H \interior \Omega = dH.
$$
We generalize this relation as
$$
X^1_H \interior \Omega_1 = \sqrt{|h_{11}|} dH,
$$
using the distinguished objects
$$
X^1_H, \; \Omega_1, H
$$
and the manifold $J^1(T,M)$. For another point of view, see also [11].

{\bf Theorem}. {\it The DEs system
$$
{d^2x^i \ov dt^2} - H^1_{11} {dx^i \ov dt} + G^i_{jk} {dx^j \ov dt}
{dx^k \ov dt} = g^{ih} g_{kj} (\na_h X^k)X^j
$$
transfers in $J^1(T,M)$ as a Hamilton DEs system with respect to the Hamiltonian
$$
H = {1\ov 2} h^{11} g_{ij} y^i y^j - f
$$
and the non-degenerate distinguished symplectic relative 2-form}
$$
\Omega = \Omega_1 \otimes dt^1, \quad \Omega_1 = g_{ij} dx^i \wedge \delta y^j
\sqrt{|h_{11}|} .
$$

{\bf Proof}. Let
$$
\te = \te_1 \otimes dt^1, \quad \te_1 = g_{ij} y^i dx^j \sqrt{|h_{11}|}
$$
be the distinguished Liouville relative 1-form on $J^1(T,M)$. We find
$$
\Omega_1 = -d\te_1.
$$
We introduce
$$
X_H = X^1_H {\delta \ov \delta t}, \; X^1_H = u^{1l} {\de \ov \de x^l} +
{\de u^{1l} \ov dt} {\pa \ov \pa y^l}
$$
as the distinguished Hamiltonian object associated to the function $H$.

The relation
$$
X^1_H \interior \Omega_1 = \sqrt{|h_{11}|} dH,
$$
where
$$
dH = h^{11} g_{ij} y^j \de y^i - h^{11} g_{ij} (DX^i) X^j dt -
h^{11} g_{ij} X^j \na_k X^i dx^k,
$$
implies
$$
g_{ij} u^{1i} \de y^j - g_{ij} {\de u^{1j} \ov dt} dx^i = dH.
$$

Consequently, it appears the PDEs system of Hamilton type
$$
\left\{ \begin{array}{l}
u^{1i} = h^{11} y^i \\ \noa
\di{\de u^{1i} \ov dt} = g^{hi} h^{11} g_{jk} X^j (\na_h X^k) \end{array}
\right.
$$
together the condition
$$
h^{11} g_{ij} (DX^i)X^j = 0.
$$

{\bf Theorem}. {\it The DEs system
$$
{d^2x^i \ov dt^2} - H^1_{11} {dx^i \ov dt} + G^i_{jk} {dx^j \ov dt}{dx^k \ov dt}
= g^{ih} g_{kj} (\na_h X^k)X^j + F_j{}^i {dx^j \ov dt} + DX^i
$$
transfers in $J^1(T,M)$ as a Hamilton DEs system with respect to the Hamiltonian
$$
H = {1 \ov 2} h^{11} g_{ij} y^i y^j - f
$$
and the non-degenerate distinguished symplectic relative 2-form
$$
\Omega = \Omega_1 \otimes dt, \; \Omega_1 = (g_{ij} dx^i \wedge \de y^j +
\omega_{ij} dx^i \wedge dx^j + g_{ij} (DX^i) dt \wedge dx^j) \sqrt{|h_{11}|},
$$
where}
$$
\omega_{ji} = g_{hi} F_j{}^h.
$$

{\bf Proof}. Let
$$
\te = \te_1 \otimes dt^1, \; \te_1 = (g_{ij} y^i dx^j - g_{ij} X^i dx^j)
\sqrt{|h_{11}|}
$$
be the distinguished Liouville relative 1-form on $J^1(R,M)$. We find
$$
\Omega_1 = - d\te_1.
$$
We denote
$$
X_H = X^1_H {\de \ov \de t}, \; X^1_H = h^{11} {\de \ov \de t} +
u^{1l}{\de \ov \de x^l} + {\de u^{1l} \ov dt}{\pa \ov \pa y^l}
$$
the distinguished Hamiltonian object of the function $H$. The relation
$$
X^1_H \interior \Omega_1 = \sqrt{|h_{11}|} dH
$$
can be written
$$
g_{ij} u^{1i} \de y^j - g_{ij} {\de u^{1j} \ov dt} dx^i + 2\omega_{ij}
u^{1i} dx^j - g_{ij} (DX^i) u^{1j} dt + h^{11} g_{ij} (DX^i)dx^j= dH,
$$
where
$$
dH = -h^{11} g_{ij} (DX^i) X^j dt + h^{11} g_{ij} y^j \de y^i - h^{11} g_{ij}
X^j (\nabla_k X^i)dx^k.
$$
Via these relations we identify a PDEs system of Hamilton type,
$$
\left\{ \begin{array}{l}
u^{1i} = h^{11} y^i \\ \noa
\di{\de u^{1i} \ov dt} = g^{hi} h^{11} g_{jk} X^j (\na_h X^k)
+ 2g^{hi} \omega_{jh} u^{1j} + h^{11} DX^i \end{array}
\right.
$$
together the condition
$$
g_{ij} (DX^i) (u^{1j} - h^{11} X^j) = 0.
$$

\section{Solutions of PDEs as Potential Maps}

All manifolds and maps are $C^\ty$, unless otherwise stated.

Let $(T,h)$ and $(M,g)$ be semi-Riemann manifolds of dimensions $p$ and $n$.
Hereafter we shall assume that the manifold $T$ is oriented. Greek (Latin)
letters will be used for indexing the components of geometrical objects
attached to the manifold $T$ (manifold $M$).

Local coordinates will be written
$$
t = (t^\al), \quad \al = 1, \ld, p
$$
$$
x = (x^i), \quad i = 1, \ld, n,
$$
and the components of the corresponding metric tensor and Christoffel symbols
will be denoted by $h_{\al \beta}, g_{ij}$, $H^\al_{\beta \gamma}$, $G^i_{jk}$.
Indices of tensors or distinguished tensors will be rised and lowered in the
usual fashion.

Let $C^\ty (T,M) = \{ \va: T \to M |\: \va \; \hbox{of class } \; C^\ty\}$.
For any $\va, \psi \in C^\ty (T,M)$ we define the equivalence relation $\va \sim \psi$
at $(t_0, x_0) \in T \ti M$, by
$$
x^i (t_0) = y^i (t_0) = x^i_0, \; {\pa x^i \ov \pa t^\al} (t_0) =
{\pa y^i \ov \pa t^\al} (t_0).
$$
Using the factorization
$$
J^1_{(t_0,x_0)} (T,M) = C^\ty (T,M)/_\sim,
$$
we introduce the jet bundle of order one
$$
J^1 (T,M) = \bigcup_{(t_0, x_0) \in T\ti M} J^1_{t_0,x_0} (T,M).
$$
Denoting by $[\va]_{(t_0,x_0)}$ the equivalence class of the map $\va$, we
define the projection
$$
\pi: J^1 (T,M) \to T \ti M, \; \pi [\va]_{(t_0,x_0)} = (t_0, \va (t_0)).
$$
Suppose that the base $T\ti M$ is covered by a systems of coordinate neighborhood
$(U \ti V, t^\al, x^i)$. Then we can define the diffeomorphism
$$
F_{U\ti V}: \pi^{-1} (U \ti V) \to U \ti V \ti R^{pn}
$$
$$
F_{UV} [\va]_{(t_0,x_0)} = \left( t^\al_0, x^i_0, {\pa x^i \ov \pa t^\al}
(t_0)\right).
$$

Consequently $J^1(T,M)$ is a differentiable manifold of dimension $p+n+pn$.
The coordinates on $\pi^{-1} (U\ti V) \subset J^1 (T,M)$ will be
$$
(t^\al, x^i, x^i_\al),
$$
where
$$
t^\al \left([\va]_{(t_0,x_0)} \right) = t^\al (t_0),
x^i \left([\va]_{(t_0,x_0)} \right) = x^i(x_0),
x^i_\al \left([\va]_{(t_0,x_0)} \right) = {\pa x^i \ov \pa t^\al} (t_0).
$$

A local changing of coordinates $(t^\al, x^i, x^i_\al) \to
(\bar t^\al, \bar x^i, \bar x^i_\al)$ is given by
$$
\bar t^\al = \bar t^\al (t^\beta), \; \bar x^i = \bar x^i(x^j), \;
\bar x^i_\al = {\pa \bar x^i \ov \pa x^j} {\pa t^\beta \ov \pa \bar t^\al}
x^j_\beta, \leqno (7)
$$
where
$$
\det \left( {\pa \bar t^\al \ov \pa t^\beta}\right) > 0, \quad
\det \left( {\pa \bar x^i \ov \pa x^j}\right) \ne 0.
$$
The expression of the Jacobian matrix of the local diffeomorphism $(7)$ shows
that the jet bundle of order one $J^1 (T,M)$ is always orientable.

Let $H^\al_{\beta\gamma}, G^i_{jk}$ be the components of the connections
induced by $h$ and $g$ respectively.
If $(t^\al, x^i, x^i_\al)$ are the coordinates of a point in $J^1(T,M)$, then
$$
x^i_{\al\beta} = {\pa^2 x^i \ov \pa t^\al \pa t^\beta} - H^\gamma_{\al\beta}
x^i_{\gamma} + G^i_{jk} x^j_\al x^k_\beta
$$
are the components of a distinguished tensor on $T \ti M$. Also
$$
\left( {\delta \ov \delta t^\al} = {\pa \ov \pa t^\al} + H^\gamma_{\al\beta}
x^i_\gamma {\pa \ov \pa x^i_\beta}, \;
{\delta \ov \delta x^i} = {\pa \ov \pa x^i} - G^h_{ik} x^k_\al
{\pa \ov \pa x^h_\al}, \; {\pa \ov \pa x^i_\al}\right),
$$
$$
\left( dt^\beta, dx^j, \; \delta x^j_\beta = dx^j_\beta - H^\gamma_{\beta\la} x^j_\gamma
dt^\la + G^j_{hk} x^h_\beta dx^k \right)
$$
are dual frames on $J^1 (T,M)$, i.e.,
$$
dt^\beta\left({\delta \ov \delta t^\al}\right) = \delta^\beta_\al, \quad
dt^\beta \left({\delta \ov \delta x^i}\right) = 0, \quad
dt^\beta \left({\pa \ov \pa x^i_\al}\right) = 0
$$
$$
dx^j\left({\delta \ov \delta t^\al}\right) = 0, \quad
dx^j  \left({\delta \ov \delta x^i}\right) = \delta^j_i, \quad
dx^j \left({\pa \ov \pa x^i_\al}\right) = 0
$$
$$
\delta x^j_\beta\left({\delta \ov \delta t^\al}\right) = 0, \quad
\delta x^j_\beta \left({\delta \ov \delta x^i}\right) = 0, \quad
\delta x^j_\beta \left({\pa \ov \pa x^i_\al}\right) = \delta^j_i \delta^\al_\beta.
$$

Using these frames, we define on $J^1(T,M)$ the induced Sasaki-like metric
$$
S_1 = h_{\al\beta} dt^\al \otimes dt^\beta + g_{ij} dx^i \otimes dx^j +
h^{\al\beta} g_{ij} \delta x^i_\al \otimes \delta x^j_\beta.
$$
The geometry of the manifold $J^1(T,M)$ was developed recently in [4].

The Lorentz world-force law formulated usually for particles [5] can be
generalized as follows:

{\bf Definition}. Let $F_\al =(F_j{}^i{}_\al)$ and $U_{\al\beta} =
(U^i_{\al\beta})$ be $C^\ty$ distinguished tensors on $T\times M$,
where $\omega_{ji\al} = g_{hi}F_j{}^h{}_\al$ is skew-symmetric
with respect to $j$ and $i$. Let $c(t,x)$
be a $C^\ty$ real function on $T \times M$. A $C^\ty$ map $\va: T \to M$ obeys the
{\it Lorentz-Udri\c ste World-Force Law} with respect to $F_\al$,
$U_{\al\beta}$, $c$ iff
$$
h^{\al\beta}x^i_{\al\beta} = g^{ij} {\pa c \ov \pa x^j} + h^{\al\beta} F_j{}^i{}_\al
x^j_\beta + h^{\al\beta} U^i_{\al\beta},
$$
i.e., iff it is a potential map of a suitable geometrical structure.

Let us show that the solutions of a system of PDEs of order one are potential
maps in a suitable geometrical structure of the jet bundle of order one.
For that we remark that any $C^\ty$ distinguished tensor field
$X^i_\al (t,x)$ on $T \times M$ defines a family of
$p$-dimensional sheets as solutions of the PDEs system of order one
$$
x^i_\al = X^i_\al (t, x(t)), \leqno (8)
$$
if the complete integrability conditions
$$
{\pa X^i_\al \ov \pa t^\beta} + {\pa X^i_\al \ov \pa x^j} X^j_\beta =
{\pa X^i_\beta \ov \pa t^\al} + {\pa X^i_\beta \ov \pa x^j} X^j_\al
$$
are satisfied.

To any distinguished tensor field $X^i_\al(t,x)$ and semi-Riemann metrics
$h$ and $g$ we associate the {\it potential energy}
$$
f: T \times M \to R, \quad f = {1 \ov 2} h^{\al\beta} g_{ij} X^i_\al
X^j_\beta.
$$
The distinguished tensor field $X^i_\al$ (family of $p$-dimensional sheets)
on $(T \times M, h + g)$ is called:

1) {\it timelike}, if $f < 0$;

2) {\it nonspacelike or causal}, if $f \le 0$;

3) {\it null or lightlike}, if $f=0$;

4) {\it spacelike}, if $f>0$.

Let $\cal E = \{x_0 \in M |\: X^i_{\al} (t,x_0) = 0$, $\forall t \in T\}$ be
the set of critical points of the system (8). If $f =$ constant, upon
rescaling on $T \ti (M \setminus \cal E)$, it may be supposed that
$f \in \{-1,0,1\}$.

The derivative along a solution of (8),
$$
{\delta \ov \pa t^\beta} x^i_\al = x^i_{\al\beta} =
{\pa^2 x^i \ov \pa t^\al \pa t^\beta} -
H^\gamma_{\al\beta} x^i_\gamma + G^i_{jk} x^j_\al x^k_\beta,
$$
produce the prolongation (system of PDEs of order two)
$$
x^i_{\al\beta} = D_\beta X^i_\al + (\na_j X^i_\al) x^j_\beta. \leqno (9)
$$
which can be converted into the prolongation
$$
h^{\al\beta} x^i_{\al\beta} = g^{ih} h^{\al\beta} g_{ij} (\na_h X^k_\al)X^j_\beta
+ h^{\al\beta} F_j{}^i{}_\al x^j_\beta + h^{\al\beta} D_\beta X^i_\al,
\leqno (10)
$$
where
$$
F_j{}^i{}_\al = \na_j X^i_\al - g^{ih} g_{kj} \na_h X^k_\al
$$
is the external distinguished tensor field which characterizes the {\it helicity}
of the distinguished tensor field $X^i_\al$.

{\bf Theorem}. {\it Any solution of PDEs system (8) is a solution of the PDEs
system (10).}

The first term in the second hand member of the PDEs system (10) is
$(grad\: f)^i$. Therefore, choosing the metrics $h$ and $g$
such that $f \in \{-1,0,1\}$, the system (10) reduces to
$$
h^{\al\beta} x^i_{\al\beta} = g^{ih} F_j{}^i{}_\al x^j_\beta +
h^{\al\beta} D_\beta X^i_\al. \leqno (10')
$$

{\bf Theorem}. {\it The solutions of PDEs system (10) are the extremals of the
Lagrangian}
$$
\begin{array}{lcl}
L &=& \di{1 \ov 2} h^{\al\beta} g_{ij} (x^i_\al - X^i_\al)(x^j_\beta - X^j_\beta)
\sqrt{|h|} = \\ \noa
&=& \left( \di{1 \ov 2} h^{\al\beta} g_{ij} x^i_\al x^j_\beta -
h^{\al\beta} g_{ij} x^i_\al X^j_\beta + f\right) \sqrt{|h|}. \end{array}
$$

{\it If $F^i_{j\al} = 0$, then this Lagrangian can be replaced by}
$$
L = \left( {1 \ov 2} h^{\al\beta} g_{ij} x^i_\al x^j_\beta + f \right)
\sqrt{|h|}.
$$

2) {\it Both Lagrangians produce the same Hamiltonian}
$$
H = \left( {1 \ov 2} h^{\al\beta} g_{ij} x^i_\al x^j_\beta - f \right)
\sqrt{|h|}.
$$

{\bf Theorem (Lorentz-Udri\c ste World-Force Law)}. {\it Every solution of
the PDEs system (8) is a horizontal potential map of the semi-Riemann-Lagrange
manifold}
$$
(T\ti M, h+g, \; N(^i_\al)_j =
G^i_{jk} x^k_\al - F_j{}^i{}_\al, \; M(^i_\al)_\beta =
- H^\gamma_{\al\beta} x^i_\gamma).
$$

{\bf Corollary}. {\it Every PDE generates a Lagrangian of order one via the
associated first order PDEs system and suitable metrics on the manifold
of independent variables and on the manifold of functions. In this sense the
solutions of the initial PDE are potential maps produced by a suitable
Lagrangian.}

{\bf Proof}. Let
$$
{\pa^rx \ov \pa (t^p)^r} = F(t^\al, x, \bar x^{(r)})
$$
be a PDE of order $r$, where $\bar x^{(r)}$ represent the partial derivatives
of $x$ with respect to $t^\al$, till the order $r$ inclusively, excepting the
partial derivative $\di{\pa^r x \ov \pa (t^p)^r}$. This equation is equivalent
to a system (8).

For the sake of simplicity, we take $r=2$. We denote $\di{\pa x\ov \pa t^\al}
= x_\al = u^\al$ and we find the partial derivatives of the functions
$(x, u^\al)$ using the system
$$
\left\{ \begin{array}{l}
x_\al = u^\al \\ \noa
u^\al_\beta = u^\beta_\al, \; \al \ne \beta \\ \noa
u^2_2 = f(t^\al, x, u^\la_\mu), \; \hbox{excepting}\; \la =\mu =2.
\end{array} \right.
$$
We shall find a PDEs system of order one with $p(1+p)$ equations, which is
of type (8). Therefore, the preceding theory applies.

\section{Covariant Hamilton Equations}

Recall that on a symplectic manifold $(Q,\Omega)$ of even dimension $q$, the
Hamiltonian vector field $X_H$ of a function $H \in \cal F(Q)$ is defined by
$$
X_H  \interior \Omega = dH.
$$
This relation can be generalized as
$$
X^\al_H \interior \Omega_\al = \sqrt{|h|}dH,
$$
using the distinguished objects $X_H, \Omega, H$ on $J^1 (T,M)$. For another
point of view, see also [11].

{\bf Theorem}. {\it The PDEs system
$$
h^{\al\beta} x^i_{\al\beta} = g^{ih} h^{\al\beta} g_{jk} (\na_h X^j_\al)
X^k_\beta
$$
transfers in $J^1(T,M)$ as a covariant Hamilton PDEs system with respect to
the Hamiltonian
$$
H = {1 \ov 2} h^{\al\beta} g_{ij} x^i_\al x^j_\beta - f
$$
and the non-degenerate distinguished polysymplectic relative 2--form}
$$
\Omega = \Omega_\al \otimes dt^\al, \quad \Omega_\al = g_{ij} dx^i \wedge
\delta x^j_\al \sqrt{|h|}.
$$

{\bf Proof}. Let
$$
\theta = \theta_\al \otimes dt^\al, \quad \theta_\al = g_{ij} x^i_\al dx^j
\sqrt{|h|}
$$
be the distinguished Liouville relative 1--form on $J^1 (T,M)$. It follows
$$
\Omega_\al = -d\theta_\al.
$$
We denote by
$$
X_H = X^\beta_H {\delta \ov \delta t^\beta}, \quad
X^\beta_H = u^{\beta l} {\delta \ov \delta x^l} + {\delta u^{\beta l} \ov
\pa t^\al}{\pa \ov \pa x^l_\al}
$$
the distinguished Hamiltonian object of the function $H$. Imposing
$$
X^\al_H \interior \Omega_\al = \sqrt{|h|} dH,
$$
where
$$
dH = h^{\al\beta} g_{ij} x^j_\beta \delta x^i_\al - h^{\al\beta} g_{ij}
(D_\gamma X^i_\al) X^j_\beta dt^\gamma - h^{\al\beta}g_{ij} X^j_\beta
\na_k X^i_\al dx^k
$$
we find
$$
g_{ij} u^{\al i} \delta x^j_\al - g_{ij} {\delta u^{\al j}\ov \pa t^\al}
dx^i = dH.
$$
Consequently, it appears the Hamilton PDEs system
$$
\left\{ \begin{array}{l}
u^{\al i} = h^{\al \beta} x^i_\beta \\ \noa
\di{\delta u^{\al i} \ov \pa t^\al} = g^{hi} h^{\al\beta} g_{jk}
X^j_\beta (\na_h X^k_\al) \end{array} \right.
$$
together the condition
$$
h^{\al\beta}g_{ij} (D_\gamma X^i_\al) X^j_\beta = 0.
$$

{\bf Theorem}. {\it The PDEs system
$$
h^{\al\beta} x^i_{\al\beta} = g^{ih} h^{\al\beta} g_{kj} (\na_h X^k_\al)
X^j_\beta + h^{\al\beta} F_j{}^i{}_\al x^j_\beta + h^{\alpha\beta}
D_\beta X^i_\alpha
$$
transfers in $J^1(T,M)$ as a covariant Hamilton PDEs system with respect to
the Hamiltonian
$$
H = {1 \ov 2} h^{\al\beta} g_{ij} x^i_\al x^j_\beta - f
$$
and the non-degenerate distinguished polysymplectic relative 2--form}
$$
\Omega = \Omega_\al \otimes dt^\al, \quad \Omega_\al =
(g_{ij} dx^i \wedge \delta x^j_\al + \omega_{ij\al} dx^i \wedge dx^j
+ g_{ij} (D_\beta X^i_\alpha) dt^\beta \wedge dx^j) \sqrt{|h|}.
$$

{\bf Proof}. Let
$$
\theta = \theta_\al \otimes dt^\al, \quad \theta_\al = (g_{ij} x^i_\al dx^j
- g_{ij} X^i_\al dx^j) \sqrt{|h|}
$$
be the distinguished Liouville relative 1--form on $J^1(T,M)$. It follows
$$
\Omega_\al = -d\theta_\al.
$$
We denote by
$$
X_H = X^\beta_H {\delta \ov \delta t^\beta}, \quad
X^\beta_H = h^{\beta\gamma} {\delta \ov \delta t^\gamma} +
u^{\beta l} {\delta \ov \delta x^l} + {\delta u^{\beta l} \ov
\pa t^\al}{\pa \ov \pa x^l_\al}
$$
the distinguished Hamiltonian object of the function $H$. Imposing
$$
X^\al_H \interior \Omega_\al = \sqrt{|h|} dH,
$$
where
$$
dH = -h^{\al\beta} g_{ij} (D_\gamma X^i_\al) X^j_\beta dt^\gamma +
h^{\al\beta} g_{ij} x^j_\beta \delta x^i_\al - h^{\al\beta} g_{ij}
X^j_\beta (\na_k X^i_\al) dx^k,
$$
we find
$$
g_{ij} u^{\al i} \delta x^j_\al - g_{ij} {\delta u^{\al j}\ov \pa t^\al}
dx^i + 2\omega_{ij\al} u^{\al i} dx^j - g_{ij} (D_\beta X^i_\alpha) u^{\al j}
dt^\beta + h^{\al\beta} g_{ij} (D_\beta X^i_\al) dx^j = dH.
$$
Consequently, we obtain the Hamilton PDEs system
$$
\left\{ \begin{array}{l}
u^{\al i} = h^{\al \beta} x^i_\beta \\ \noa
\di{\delta u^{\al i} \ov \pa t^\al} = g^{hi} h^{\al\beta} g_{jk} X^j_\beta
(\na_h X^k_\al) + 2g^{hi} \omega_{j h \al} u^{\al j}
+ h^{\alpha\beta} D_\beta X^i_\alpha \end{array} \right.
$$
together the condition
$$
g_{ij} (D_\gamma X^i_\al) (u^{\al j} - h^{\al\beta} X^j_\beta) = 0.
$$

\end{document}